\documentclass{amsart}

\newcommand{\RR}{\mathbf R}
\newcommand{\CC}{\mathbf C}
\newcommand{\HH}{\mathbf H}
\newcommand{\rep}{\mathbf}
\title{Integer versions of Yang--Mills theory}
\author{Robert Arnott Wilson}

\date{First draft: 16th February 2022. This version: 8th March 2022.}
\address{Queen Mary University of London}
\email{r.a.wilson@qmul.ac.uk}

\begin{document}
\begin{abstract}
In a recent paper, 't Hooft asks for an integer version of Yang--Mills theory, in the belief that
this is the way the universe really 
is at the Planck scale. Specifically, he asks for an integer version of the 
gauge group of the standard model. Such groups were completely classified more than 100 years ago,  
and here I 
detail the group theory and representation theory (and some 
of the physics) of the specific case that I believe answers
't Hooft's question.
\end{abstract}

\maketitle

\section{Introduction}
\subsection{Finite groups}
Yang--Mills theory \cite{YangMills} is the cornerstone of the standard model of particle physics, 
and its immense success poses an enormous challenge to any new theories that attempt to go
beyond the standard model in any way. In a recent paper \cite{tHooft22} dedicated to Yang, 't Hooft suggests that 
on the Planck scale reality should be discrete, and that a discrete (and therefore finite)
version of $U(1)\times SU(2)\times SU(3)$ will be required in order to describe the fundamental theory.
He then says he has no idea which finite group to use.
On the other hand, these finite groups have been well-known to mathematicians for more than a hundred years \cite{Blichfeldt},
and 
over a period of time
I have examined all 
of them closely in the hope of being able to answer
't Hooft's question.
A number of possibilities are discussed in \cite{finite,octahedral,gl23model,icosa} but only the last of these cases 
survives this detailed scrutiny. 

Once it is clear which group we have to use, it is much easier to build the appropriate model,
with much less speculation, and much more confidence that this is the one that must work, if
any of them does. The purpose of this paper, therefore, is to go back to the beginning,
and explain from first principles what the group is, where it comes from, and what 
the implications are.

\subsection{Clifford algebras}
The concept of spin is central to much of algebra, geometry and quantum theory,
and is constructed mathematically using Clifford algebras \cite{Porteous}.
There has always been a dilemma as to whether the signature of spacetime 
is $(3,1)$ or $(1,3)$.
Naively, 
it seems it should not matter, but unfortunately it does, because the two Clifford algebras are
quite different. This is even true in the non-relativistic case, where the signatures $(3,0)$ and $(0,3)$ also give
non-isomorphic Clifford algebras:
\begin{align}
Cl(3,0) \cong M_2(\CC),\quad&
Cl(0,3)\cong \HH+\HH,\cr
Cl(3,1) \cong M_4(\RR),\quad&
Cl(1,3) \cong M_2(\HH).
\end{align}
The standard approach to this problem is to ignore it, and work instead with the complex matrix algebra
$M_4(\CC)$, known as the Dirac algebra. 
Since $M_4(\CC)$ contains both $M_4(\RR)$ and
$M_2(\HH)$, it then becomes unnecessary to make a decision.

This procedure, however, only evades the problem, and does not solve it. There is a real 
question here of 
which of these Clifford algebras describes the 
transformation properties of  spins of elementary particles in the real world. 
The answer depends 
on 
properties of neutrinos, such as whether they are Dirac or Majorana particles, so that it is not easy to give a definitive answer
on experimental grounds. Nevertheless, the principle of unitarity would seem to imply that we require $M_2(\HH)$
to act on the Dirac spinors.
Of course, this structure only allows us to implement half of the Dirac algebra, so it may not necessarily contain
enough information to cover everything that is done with the Dirac algebra in the standard model.
In particular, we have to choose between $\gamma_0$ and $\gamma_5$, because we cannot have both at the same time.
Whichever  
we choose, its eigenvalues are $1$ and $-1$, which means that it acts as 
a quaternionic reflection on the space $\HH^2$ of Dirac spinors. It may therefore be useful to study the two-dimensional quaternionic
reflection groups.

\section{Finite group theory and fundamental particles}
\subsection{Quaternionic reflection groups}
The Dirac matrices themselves generate a quaternionic reflection group, of order $32$, in which there are ten reflections 
\begin{align}
\pm\gamma_0,\quad \pm\gamma_0\gamma_1, \quad
\pm\gamma_0\gamma_2, \quad\pm\gamma_0\gamma_3, \quad\pm\gamma_1\gamma_2\gamma_3.
\end{align}
However, this group does not show any asymmetry between the two quaternion coordinates, so does not capture enough of the
known structure of physical spin for our purposes.
A complete classification of finite quaternionic reflection groups was obtained by Cohen \cite{Cohen}. The interesting cases are listed
 in \cite[Table III]{Cohen}, where the $2$-dimensional examples 
 are denoted $O_1$, $O_2$, $O_3$, $P_1$, $P_2$ and $P_3$. 
 
 Those of type $P$ all contain the above
group of order $32$ generated by the Dirac matrices $\gamma_\mu$, extended by a dihedral, alternating or symmetric group on $5$ points.
They were described in detail in \cite{icosians}, and the connection to Dirac matrices is pointed out in \cite{icosa}, but they also have
a symmetry between left-handed and right-handed spinors.
Those of type $O_2$ and $O_3$ were mentioned very speculatively in \cite{finite} as possibly having a connection to 
quantum theory, but again there is
a left-right symmetry.
Thus the only one of the six that exhibits an asymmetry between left-handed and right-handed spinors is $O_1$. It therefore seems
likely that this is the one that will be most useful for applications in quantum theory.

As an abstract group, $O_1$ is isomorphic to the binary icosahedral group, of order $120$. Possible uses of this group in quantum
physics have been
discussed at length in \cite{icosa}, but the quaternionic reflection group property renders much of that speculation redundant, and provides a
much more direct route from the mathematics to the 
applications. In particular, there is a distinguished set of $120$ Dirac spinors, namely the `roots'
which define the reflections, as well as a distinguished set of $120$ matrices in the Dirac algebra. This enables us to 
define sets of elementary fermions and bosons, and to investigate the relationships between them in detail. Of course, it is not obvious 
at this stage that
the elementary particles defined by this mathematical model actually match up to the elementary particles in the real world.

In order to define the required objects, we first need 
two quaternions $\omega$ and $\phi$ that satisfy the relations
\begin{align}
\phi^2=(\omega\phi)^2=\omega+\omega^2=-1. 
\end{align}
Then the group of scalars generated by $\omega$ and $\phi$ has order $12$, and consists of all elements $\omega^a\phi^b$, where we can suppose
that $a=0,1,2$ and $b=0,1,2,3$. The $120$ roots can then be taken as these $12$ scalar multiples of the following $10$ spinors of 
(squared) norm $3$, where again we can assume $c=0,1,2$:
\begin{align}
(\omega-\omega^2,0), \quad ((\phi+1)\omega^c,1), \quad (\omega^c,(\phi-1)\omega), \quad (\omega^c,(\phi-1)\omega^2).
\end{align}
Each root $\mathbf r$ then defines a reflection $R_\mathbf r$ via the formula
\begin{align}
R_\mathbf r : \mathbf x \mapsto \mathbf x - \frac{\mathbf x.\mathbf r}{\mathbf r.\mathbf r}(1-\omega)\mathbf r.
\end{align}
Notice that this formula does not change when we replace $\mathbf r$ by $\pm\omega^a\mathbf r$, but if we replace $\mathbf r$ by
$\pm\phi\mathbf r$, then the factor $(1-\omega)$ becomes $-\phi(1-\omega)\phi = (1-\omega^2)$. 
Hence we obtain $20$ reflections of order $3$, coming in
$10$ inverse pairs. The remaining $100$ elements of the group can all be written as products of two reflections.

In matrix notation, we can take generators 
\begin{align}
& f:= \frac{\omega-\omega^2}{3}\begin{pmatrix}1 & \omega^2\phi-\omega\cr \omega^2\phi-\omega^2 & \omega\phi\end{pmatrix},\quad\cr
&g:=\begin{pmatrix}\omega&0\cr0&1\end{pmatrix}, \quad h:=\begin{pmatrix}\phi&0\cr0&\phi\end{pmatrix}
\end{align}
such that the subgroup of diagonal matrices is a maximal subgroup, of order $12$, generated by $g$ and $h$.
These matrices satisfy the relations
\begin{align}
&f^2=(gh)^2=h^2=-1, \cr &g^3=(fg)^3=(fh)^3=1.
\end{align}

\subsection{Classification of particles}
Of particular interest for 
potential applications is the fact that $g$ acts only on the first quaternion coordinate, which might therefore correspond 
in some way to the
`left-handed' spinor in the standard model. Since there is no corresponding symmetry acting on the right-handed spinor, this is
a natural way in which the mathematics could model the `chirality' of the weak interaction. In this way we can distinguish the (Dirac)
neutrinos as the $12$ multiples of the left-handed spinor
$(\omega-\omega^2,0)$. This includes a factor of $3$ for the three generations, and a factor of $4$
which is sufficient to distinguish both neutrinos from antineutrinos, and spin up from spin down. 

The other nine types of particles also come in $12$ states each, and can presumably be identified with the three generations of electrons,
and six flavours of quarks, in some way. This gives us $72$ quark states altogether, which is the same number as in the standard model, if we
distinguish six flavours, three colours and three anti-colours, and spin up/down. Each can also be split into left-handed and right-handed
parts, but these are not particle states in the same sense. Finally we have $36$ electron states, which is three times as many as in the
standard model. Hence this model can distinguish six spin states for each electron, compared to the simple pair of spin up/down
in the standard model.

In other words, there is a `hidden variable' that takes one of three values, and that is not in the standard models of quantum mechanics
and particle physics. 
This does not 
contradict Bell's Theorem \cite{Bell,Rae}, which only forbids \emph{continuous} local hidden variables.
Two entangled electrons can therefore share this 
hidden quantum number, and use it to determine the binary spin state
as measured in any particular experiment. This is extra information that is available to the electrons, but is unknown to quantum mechanics,
and may be sufficient to obviate the need for any hypothetical superluminal transfer of information between entangled particles.

Presumably this hidden  
quantum number can also be carried by photons, in order to explain the properties of
entangled photon polarisations, but for that we need a full classification of the bosonic fields.
For bosons, then, we need to map from $\omega$ and $\phi$ in $SU(2)$ into the orthogonal group $SO(3)$, where we obtain a
group of order $6$, that is the rotation symmetry group of a triangular antiprism. In other words, a photon in this model has six intrinsic 
polarisation states rather than just two, again acquiring just enough additional information to explain the observed properties of
entangled photons without superluminal transfer of information.

\subsection{Quantum fields}
Each of the fundamental fermions generates a reflection, that can be written as a $2\times 2$ quaternion matrix, embedded in the
Dirac algebra, and therefore has an interpretation as a quantum field \cite{WoitQFT}
of some kind. There are $20$ such matrices, and together they
generate a $16$-dimensional algebra. The neutrinos, for example, generate a $3$-dimensional algebra with basis $1,g,g^2$. 
This is a commutative algebra, isomorphic to $\RR+\CC$ with the following generators:
\begin{align}
&(1+g+g^2)/3 = \begin{pmatrix}0&0\cr0&1\end{pmatrix},\quad\cr
&(2-g-g^2)/3 = \begin{pmatrix}1&0\cr0&0\end{pmatrix},\quad
g-g^2=\begin{pmatrix}\theta&0\cr0&0\end{pmatrix},
\end{align}
where $\theta=\omega-\omega^2$ is pure imaginary, with $\theta^2=-3$. We can also adjoin $h$ to this algebra, in order to
incorporate the spin of the neutrinos, and hence obtain a $6$-dimensional algebra isomorphic to $\CC+\HH$.

To see this isomorphism, adjoin the generators
\begin{align}
&h(1+g+g^2)/3 = \begin{pmatrix}0&0\cr0&\phi\end{pmatrix},\cr
&h(2-g-g^2)/3 = \begin{pmatrix}\phi&0\cr0&0\end{pmatrix},\quad
h(g-g^2)=\begin{pmatrix}\phi\theta&0\cr0&0\end{pmatrix}.
\end{align}
The corresponding Lie algebra \cite{Zee} is
\begin{align}
gl(1,\CC)+gl(1,\HH) = gl(1,\RR) + u(1) + gl(1,\RR)+su(2),
\end{align}
which we would like to interpret in relation to the gauge group $SU(2)_L\times U(1)_Y$. If this is a viable interpretation,
then $u(1)_Y$ must be generated by $h(1+g+g^2)/3$ acting on the right-handed spinor only, while $su(2)_L$ is
generated by 
\begin{align}
h(2-g-g^2)/3, \quad g-g^2, \quad h(g-g^2), 
\end{align}
acting on the left-handed spinor only. 

Certainly this would show us the symmetry-breaking of the electroweak gauge group
$SU(2)\times U(1)$ quite clearly. It also shows 
that the standard formulation involves taking quantum superpositions of
the discrete particles implied by the finite group. 
To put this another way, the discrete particles that are actually observed (the $Z$ and $W$ bosons
and the photon) do not fit neatly into the Lie algebras of the standard model. Indeed, $u(1)_{em}$ is most naturally defined
to be generated by $h$, so that $h$ itself can be identified as an abstract photon, divorced from spacetime. If so, then the
symmetries imply that $hg$ and $gh$ are also photons, which gives us the extra photon states already conjectured as
being required to explain the experimental properties of the measurement of photon polarisation. This leaves us with $1$ to
denote the $Z$ boson, and $g,g^2$ for the $W^+$ and $W^-$ bosons.

However, this interpretation of photons is not consistent with the standard model, which would require $h,hg,gh$ to be
`primordial' massless weak bosons, so that the photons must go somewhere else. It is reasonable to assume that the massless bosons
correspond to the elements of order $4$ in the group, of which there are $15$ up to sign. 
Any choice of signs gives us a basis for the adjoint representation of $SL(2,\HH)$, and the symmetry group
$\langle g,h\rangle$ should describe the division into electrodynamics and the weak and strong forces. 

With respect to $\langle g,h\rangle$ there are three orbits of lengths $3+6+6$ on the $15$ pairs $\pm x$ of elements of order $4$
in the full group:
\begin{align}
&h,gh,hg,\cr & f,f^g,f^{g^2}, f^h,f^{gh},f^{hg},\cr & f^{ghf},f^{hgf},f^{ghfg},f^{hgfg},f^{ghfg^2},f^{hgfg^2}.
\end{align}
The first orbit consists of elements inside the subgroup, while the second orbit consists of elements
in the two copies of the quaternion group $Q_8$ that are normalised by $g$.  The third orbit consists of three pairs, with the
property that the product of the two elements in each pair lies in the original group. A possible interpretation is that the second orbit
consists of photons, with conjugation by $h$ corresponding to a change of polarisation. This would force the third orbit to
consist of coloured gluons, and the first orbit of colourless gluons, extended from two to three dimensions in order to avoid the
need for quantum superposition. 

\subsection{Elements of odd order} 
Let us now look at some of the other elements of the group, that represent massive bosons. 
We have seen that $g$ and $g^2$ might represent the $W^+$ and $W^-$ bosons, which suggests that the other elements
of order $3$ might represent 
massive bosons related to the strong force and/or electromagnetism.
Up to inversion, the elements of order $3$ fall into orbits of lengths $1+3+6$ under the symmetry group $\langle g,h\rangle$.
Conjugation by $g$ should be a weak force process that acts on fermions to pair the neutrino with the electron,
the proton with the neutron, and so on, so approximately preserves the mass, while appearing to mix up the charges.
Presumably it has a similar effect on bosons.
Conjugation by $h$ on the other hand preserves the mass exactly, but appears to negate the charge on the bosons.
Similarly, inversion also appears to fix the mass and negate the charge.

This suggests the splitting $3+6$ is related to a splitting of massive strong force mediators into pions and kaons. 
The former case then consists of the six elements 
\begin{align}
fh, (fh)^g, (fh)^{g^2}, \quad hf, (hf)^g, (hf)^{g^2},
\end{align}
which gives two mutually inverse states for each of the three pions. 
This doubling up of states is required in order to avoid having to express the neutral pion as a quantum superposition of
$u\bar u$ and $d\bar d$. The states are distinguished (in theory) by which of the two colourless gluons is used to glue the
quark and anti-quark together. Since these gluons cannot be directly observed, this model gives the correct number of pions in total.
Both operations of conjugation by $h$, and inversion, act to negate the charge.
In the other case, therefore, we should expect to see no more than six experimentally distinguishable particles among the
$12$ group elements. 

These group elements are
\begin{align}
g^f, g^{fg}, g^{fg^2}, \quad g^{fh}, g^{fgh}, g^{fg^2h}
\end{align}
and their inverses. The mathematics allows either two charged and four neutral particles here, or four charged and two neutral,
but experiment certainly favours the former.
Our assumptions allow at most three distinct mass values, with four particle states each. It seems more likely that 
there are only two distinct masses,
giving two charged and four neutral kaons, as in the standard model, but without the need for quantum superposition.
There may, of course, be other possible interpretations.

The remaining elements of the group have order $5$ modulo sign, and form $6$ cyclic groups, with generators
\begin{align}
gfh,fhg,ghfg, \quad ghf,hfg,gfhg.
\end{align}
In the case of an element of $5$, however, we should not expect its mass to be related in any simple way
to the mass of its square, so that it would be reasonable to expect four distinct charged pairs and four
masses of neutral bosons. 
This means we have enough room for all the remaining charged pseudoscalar mesons, of types $D^\pm$, $D_s^\pm$, $B^\pm$ and
$B_c^\pm$, since the top quark does not hadronise.

There does not seem to be enough room for all the neutral pseudoscalar mesons, however, and perhaps the various types of eta meson
must be explained in a different way. Then we have three neutral types $D^0$, $B^0$ and $B^0_s$, which leaves room for $Z^0$
to take the fourth available mass value, so that it can have spin $1$ and the identity element of the group then has spin $0$
and can be allocated to the Higgs boson.
Alternatively, perhaps $h$ is not constrained to fix the mass of neutral particles in the way that it seems to be constrained to fix
the mass of charged particles. In that case, the neutral mesons may be distributed in quite a different way from the way I have suggested.

To summarise, in this context of a potential unification of electromagnetism with both the weak and strong forces,
as well as a potential unification of the leptons and quarks,
the price of 
separating the 
algebra $M_4(\CC)$ into two separate parts $M_2(\HH)$ and $M_4(\RR)$, with significantly different
physical interpretations, seems a price worth paying. 
This also gives us a useful opportunity to re-write the relationship between the quantum `quaternionic world' of $M_2(\HH)$
and the 
relativistic `real world' of $M_4(\RR)$.
In order to investigate this possibility further, we need to study the representation theory  \cite{JamesLiebeck}
of the binary icosahedral group
in more detail.

\section{Representation theory and effective field theory}
\subsection{Irreducible representations}
It is well-known \cite{Klein,Klein1} that the binary icosahedral group is a subgroup of $SU(2)$, and therefore of the real
quaternions \cite{DuVal}. A complete list of its elements in quaternionic notation is given in \cite{icosians}, 
reprinted in \cite{SPLAG}, from which we can
easily read off suitable generators as
\begin{align}
f\mapsto i, \quad g\mapsto (-1+i+j+k)/2, \quad h\mapsto i-\sigma j-\tau k, 
\end{align}
where $\tau=(1+\sqrt5)/2$ and $\sigma=(1-\sqrt5)/2$.
All the irreducible representations
can then be labelled by the `spin' of the corresponding representations of $SU(2)$,
but at this stage we do not know what interpretation to give to this copy of $SU(2)$, so we do not know what type of spin or isospin
this might be. 
Let us therefore call it hyperspin, and define the above representation
to have hyperspin $1/2$. 

Then the representations with hyperspin up to $7/2$ 
include all the
irreducibles, labelled here with their complex dimensions:
\begin{align}
\begin{array}{cc}
\mbox{Hyperspin} & \mbox{Decomposition}\cr\hline
0 & \rep1\cr
1 & \rep3a\cr
2 & \rep5\cr
3 & \rep3b+\rep4a
\end{array}
\qquad
\begin{array}{cc}
\mbox{Hyperspin} & \mbox{Decomposition}\cr\hline
1/2 & \rep2a\cr
3/2 & \rep4b\cr
5/2 & \rep6\cr
7/2 & \rep2b+\rep6
\end{array}
\end{align}
The representation $\rep4b$ is the quaternionic reflection representation described above, and $\rep2b$ can be obtained from
$\rep2a$ by changing the sign of $\sqrt5$. 

The other representations can be obtained from the tensor products
\begin{align}
&\rep2a\otimes\rep2a=\rep1+\rep3a\cr
&\rep2a\otimes\rep2b=\rep4a\cr
&\rep2b\otimes\rep2b=\rep1+\rep3b\cr
&\rep2b\otimes\rep3a=\rep6\cr
&\rep2b\otimes\rep4b=\rep3b+\rep5
\end{align}
The integer hyperspin representations are all real, so 
generate real matrix algebras
\begin{align}
\RR + M_3(\RR)+M_3(\RR)+M_4(\RR)+M_5(\RR),
\end{align}
while the half-integer hyperspin representations are all quaternionic, so generate 
quaternionic matrix algebras
\begin{align}
\HH+\HH+M_2(\HH)+M_3(\HH).
\end{align}
In particular, there is an obvious copy of $M_4(\RR)$
with which to attempt to reproduce the other half of the Dirac algebra,
acting on $\rep4a$. 

In order to do this, we will need to understand how the finite group acts on this matrix algebra. This action is described by
the tensor product
\begin{align}
\rep4a\otimes\rep4a=\rep1+\rep3a+\rep3b+\rep4a+\rep5.
\end{align}
Notice that this is equivalent to the Dirac algebra in the form
\begin{align}
\rep4b\otimes\rep4b=\rep1+\rep3a+\rep3b+\rep4a+\rep5.
\end{align}
In the standard model, this is a complex tensor product, and therefore a complex representation, but here it is a
quaternionic tensor product, and therefore a real representation.

This remarkable coincidence implies that at the quantum level, that is at the level of the finite symmetries of the elementary particles,
there is an exact equivalence between the $4\times 4$ real matrices and the $2\times 2$ quaternion matrices. But this equivalence
is only an equivalence of \emph{representations}, not an equivalence of \emph{algebras}. 
Hence we can add the matrices together with impunity, but we must be very careful about multiplying them together.
That is, we can \emph{either} multiply matrices in $M_2(\HH)$, as is done in the Feynman calculus in the standard model
of particle physics, \emph{or} multiply matrices in $M_4(\RR)$, as is done in general relativity, but we cannot do both at the
same time. In particular, we will need to pay particular attention to the differences between multiplication in the
subgroup $SO(3,1)$ of $M_4(\RR)$ and multiplication in the subgroup $SL(2,\CC)$ of $M_2(\HH)$.

\subsection{Gauge groups}
Gauge groups are used everywhere in physics to describe coordinate transformations in the mathematical theory that have no effect
on the physical phenomena. These include the Lorentz group $SO(3,1)$ and the group $GL(4,\RR)$ of general covariance,
and the Dirac group $SL(2,\CC)$ (often also 
called the Lorentz group), but elsewhere in particle physics only unitary 
groups are used. The underlying mathematical (as opposed to physical) 
reason for this is that representations of finite groups are always unitary. However,
there are really three types of representations of finite groups, that is real (orthogonal), complex (unitary) and quaternionic (also
called symplectic or pseudoreal). In the case of the binary icosahedral group, only orthogonal and symplectic
representations arise, and there are no unitary groups.

The real representations have orthogonal gauge groups
\begin{align}
SO(3)\times SO(3)\times SO(4)\times SO(5)
\end{align}
and the quaternionic representations have symplectic gauge groups
\begin{align}
Sp(1)\times Sp(1) \times Sp(2)\times Sp(3).
\end{align}
Since there are natural two-to-one maps from $Sp(1)\times Sp(1)$ onto $SO(4)$, and from $SO(4)$ onto $SO(3)\times SO(3)$, and from
$Sp(2)$ onto $SO(5)$, it is only necessary to use the symplectic groups. In the standard model, the quaternionic structure of $\rep4b$
that defines the three generations of fermions is not used, so that only the (non-relativistic) spin group $SU(2)$ remains. The same seems to be
true for $\rep6$, reducing the gauge group to $SU(3)$, and for one of $\rep2a$ or $\rep2b$, reducing the gauge group to $U(1)$.
In total, this reduces the gauge group from $37$ dimensions to $15$, leaving $22$ arbitrary parameters behind,
split roughly into three types as 
\begin{align}
2+(1+3+3)+(1+6+6).
\end{align}

Hence the model provides a plausible \emph{mathematical} source for the unexplained parameters of the standard model.
This is of course not the same as providing a plausible \emph{physical} explanation, which is likely to be a great deal more difficult.
At the same time, the model provides a selection of extra gauge groups with which to gauge these parameters.
As it stands, there is some ambiguity as to whether to allocate 
weak $SU(2)$ 
to a subgroup of $Sp(2)$ acting on $\rep4b$, and non-relativistic spin to
$Sp(1)$ acting on $\rep2a$, or \emph{vice versa}. It is not yet clear
which version will give a better match to the real world, so 
we need to keep an open mind on this question. 

Either way, electro-weak mixing in the standard model identifies $\rep2a+\rep2b$ with $\rep4b$,
which is incompatible with the action of any of the triplet symmetries. One of these is the generation symmetry,
so that the standard model Dirac algebra can only deal with one generation at a time. Another is the triplet
symmetry of spin states, so that the standard model Dirac algebra can only distinguish two spin states rather than six,
and therefore cannot explain the experimental properties of entangled spin states.
It seems clear, therefore, that distinguishing two completely different types of `Dirac spinor', in $\rep2a+\rep2b$ and $\rep4b$,
will be necessary in order to go beyond the standard model.
Indeed, there is a real ambiguity in the standard model as to whether the Dirac spinor should be acted on by $Spin(1,3)$,
in which case it must lie in $\rep4b$, or by $Spin(4)$, in which case it must lie in $\rep2a+\rep2b$.
Furthermore, the Dirac equation, essential for the quantum mechanical definition of mass, cannot be written
on either $\rep4b$ or $\rep2a+\rep2b$ alone, because it requires a complex scalar that commutes with the
Dirac matrices.

To make the gauge groups of the standard model, therefore, we start with the spin group $Sp(1)$ acting on a non-relativistic
spinor in $\rep2a$, and the weak $SU(2)_L$ as a subgroup of $Sp(2)$ acting on 
elementary particle labels in $\rep4b$, or possibly \emph{vice versa}. 
We also have a copy of $U(1)$ acting on $\rep2b$, and thereby also
defining complex structures on 
\begin{align}
&\rep2b\otimes \rep3a=\rep6,\cr &\rep2b\otimes \rep3b=\rep2b+\rep4b, 
\end{align}
so that
we obtain three copies of $U(1)$ with various possible interpretations. However, these three copies are really all the same, so that
although we might interpret one copy as $U(1)_Y$ acting on $\rep2b$, and therefore commuting with $SU(2)_L$ acting on $\rep4b$,
what we really have is a copy of $U(1)$ that also acts on $\rep4b$. But the latter is not guaranteed to commute with $SU(2)_L$,
and 
most likely it does not, so that it breaks the symmetry of $SU(2)_L$ and acquires a new interpretation as $U(1)_{em}$.

In order to go beyond the standard model and incorporate the triplet symmetries, it is necessary and sufficient to extend from
$U(1)$ to $SU(2)$ acting on $\rep2b$. This now provides a quaternionic structure on $\rep2b$, which then transfers to
a quaternionic structure on $\rep4b$ and $\rep6$, with which we can construct the required larger gauge groups $Sp(2)$ and $Sp(3)$.
By extending the gauge group $U(1)_Y\times SU(2)_L\times SU(3)_c$ of the standard model to
$Sp(1)\times Sp(2)\times Sp(3)$ we extend  QED and the weak and strong forces to include all three generations.
Of course, this is not in any sense a Grand Unified Theory, it is merely a suggestion for how the gauge groups of
such a theory might be derived from a few basic physical principles. Building the theory itself is a much more difficult undertaking.

\subsection{Effective field theories}
If 't Hooft \cite{tHooft17} is right, and a discrete model along the lines I have sketched
is the fundamental design of the universe at the Planck scale, then at larger scales 
much of the discrete structure appears to us to be continuous, and the continuous gauge
groups form the basis for the effective field theories that are used in practice. The simplest
way to derive the effective field theories from the discrete model would seem to be to use
$\rep4b$ for the Dirac spinor, acted on by $Spin(1,3)=SL(2,\CC)$ contained in the quaternionic 
part of the Dirac algebra, generated by $\gamma_\mu$, for $\mu=0,1,2,3$. 

The standard model then converts $Spin(1,3)$ acting on $\rep4b$ to $Spin(4)$ acting on $\rep2a+\rep2b$
by the expedient of multiplying the `time' coordinate by $i$, thereby replacing $\gamma_0$ by $\gamma_5$,
and splitting the spinor into two chiral pieces. At the same time, the triplet symmetries of the
finite group are destroyed, so that flavour and colour symmetries have to be incorporated `by hand'.
The physical nature of this chirality comes from the natural map from $Spin(4)$ acting on $\rep2a+\rep2b$ onto
$SO(4)$ acting on $\rep2a\otimes\rep2b=\rep4a$. Finally, we convert back to $SO(3,1)$ by converting again
between real and imaginary time.

In effective field theories, it is assumed that this process of getting from $Spin(1,3)$ to $SO(1,3)$ is
given by a natural mathematical two-to-one map, but if it is obtained in the way I have described then there
is nothing natural about it. Of the three necessary steps, only the second is natural, and both the other two
involve choices. In 't Hooft's terminology, these choices are made not by God, but by us. For example, in the last step,
God's computer program works with a finite subgroup of $SO(4)$, and knows nothing of $SO(3,1)$. 

In order to clarify the issues, we can separate the two choices to some extent  by supposing the identification
of $\rep2a+\rep2b$ with $\rep4b$ to have already been done, using the Dirac equation,
so that electrodynamics can be expressed in terms of 
\begin{align}
(\rep2a+\rep2b)\otimes(\rep2a+\rep2b) &= \rep1+\rep4a+\rep3a+\rep3b+\rep4a+\rep1.
\end{align}
This allows us to implement the usual interpretation of the Dirac algebra, in which $\rep4a$ is a vector
representation of $SO(3,1)$, so can be interpreted as $4$-momentum and/or spacetime position
as appropriate. The other half of the algebra consists of $\rep3a+\rep3b$ standing for the electromagnetic field,
plus two scalars for mass and charge. But in the fundamental theory, one copy of $\rep4a+\rep1$, representing
$4$-momentum and mass, is replaced by $\rep5$, as is necessary in order to allow the weak interaction
to change the mass of elementary particles without changing the total energy.

In particular, we need to choose a copy of $SO(3,1)$ in order to define mass, without which our effective
field theories do not work. But God's computer program does not have a copy of $SO(3,1)$, and therefore does not
have a concept of mass. Therefore, mass cannot be a \emph{fundamental} physical concept, but must be
emergent. In other words, mass is not a cause, but an \emph{effect}. In particular, mass cannot be the
ultimate {cause}
of gravity, however much it may seem like that to us. Perhaps this is why it has proved so difficult to
quantise gravity: have we misunderstood its cause? I emphasise that this conclusion is not as ridiculous
as it sounds: it is a necessary consequence of following 't Hooft's line of reasoning to its logical
conclusion.

\section{Consequences for mass and gravity}
\subsection{The measurement of mass}
If we accept this conclusion, then it becomes much easier to understand such mysteries as why the electron mass
is so small, and why the proton and neutron masses are so close to each other. These facts then cease to be
fundamental properties of the universe, but properties that exist only in the eye of the beholder. This on its own
does not explain the facts, but it does show that we have been looking in the wrong place for the answer. We need
to look more closely into our own eyes. If the electron/proton/neutron mass ratios have no fundamental meaning,
but only a practical meaning in effective field theory, then the near equality of proton and neutron masses can be
put down to pure coincidence and nothing more.

Moreover, the mass ratios that we use in our effective field theories depend on our choice of $SO(3,1)$,
which in practice is determined by our assumption that the laboratory frame of reference is near enough
inertial that it doesn't matter. But the laboratory frame of reference is not inertial, so that the actual copy of
$SO(3,1)$ that we use is crucially dependent on such accidents as the relative lengths of the day, the month
and the year, the angle of tilt of the Earth's axis, the eccentricity of the orbits of the Earth and the Moon,
and many other factors. That does not mean that the masses change when any of these parameters changes,
because we are free to choose a `standard' copy of $SO(3,1)$ that only approximately
describes the laboratory frame of reference, so that the practical variability can be moved into the
identification of $SU(2)\times SU(2)$ with $SL(2,\CC)$ instead. In order to investigate whether this choice of
$SO(3,1)$ actually matters or not, we need to look at the history of this choice, and see whether it has left
its imprint on the parameters of the standard model, in the form of suspicious coincidences.

To demonstrate that such coincidences do exist, consider the coincidence between the neutron/proton mass ratio
\begin{align}
m(n)/m(p) \approx 1.001378
\end{align}
and the following formula based on the average number of days in a year:
\begin{align}
1+1/2\times 365.24 \approx 1.001369.
\end{align}
Clearly this is a pure coincidence, without physical meaning. To claim otherwise would be absurd. If we look at the 
electron as well, we have
\begin{align}
m(e)/m(p) \approx .000544617
\end{align}
and considering the angle of tilt of the Earth's axis to be around $23.44^\circ$ we have the following formula:
\begin{align}
\sin(23.44^\circ)/2\times 365.24 \approx .000544558.
\end{align}
Clearly this is a pure coincidence, without physical meaning. To claim otherwise would be absurd.

In any case, the angle of tilt of the Earth's axis varies considerably, but the electron/proton mass ratio does not.
We can ask what angle makes the coincidence exact?
We need
\begin{align}
&\sin\theta \approx 2\times 365.24\times .000544617 \approx .3978318\cr
&\theta \approx 23.442704^\circ\approx 23^\circ 26' 33.7''.
\end{align}
Now we can check the astronomical almanac \cite{almanac} to find when this value of the angle of tilt
was attained. This has happened only three times since the end of the last ice age,
in August 1957, June 1963 and March 1973. Is it a coincidence that this is 
the time that the standard model of
particle physics was being built? Or is this real evidence that \emph{we} chose the 
electron/proton mass ratio, not God?

\begin{quote}
\textit{The fault, dear Brutus, lies not in the stars, but in ourselves, \ldots} 
\end{quote}

If we accept 't Hooft's argument, then God's computer program that runs the universe does not contain
the electron/proton mass ratio. But our effective field theories that describe how we see the universe
could not work without 
this mass ratio. Hence we are forced by our desire to have
effective field theories to pick a value for it. 
 But God does not care what value
we pick, so we can (effectively) pick any value we like. The fact that we (effectively) picked the current value
in the 1973 CODATA revision of fundamental physical constants \cite{1973}, and have not materially changed it since,
does not of itself give this particular value any fundamental physical meaning. 
I emphasise that this argument does \emph{not} in any sense `explain' the electron/proton mass ratio, but
only suggests that an explanation might eventually be found in a quantum theory of gravity that mixes with
the other forces along the lines indicated above.

An effective field theory can be built
with \emph{any} choice of this parameter, 
and the job of CODATA is to fix a complete consistent set of parameters that work, not to
determine a `correct' or `universal' value of any single parameter.
It is noticeable, for example, that the coincidence of the neutron/proton mass ratio cannot be made exact
by fixing a date. This implies that the process of determining this mass ratio was historically more complicated than that
of determining the electron/proton mass ratio, so that guessing a simple empirical formula does not work.

\subsection{The equivalence principle}
The above discussion strongly suggests that not all definitions of mass are equivalent. This question is usually
framed in Newtonian rather than relativistic form, that is as the question of equivalence or otherwise of
inertial mass defined by
\begin{align}
m:= F/a
\end{align}
and (active) gravitational mass, defined by
\begin{align}
M:= ar^2/G
\end{align}
on the assumption (strongly supported by near-Earth experiment) that passive gravitational mass is the
same as inertial mass. Whilst it is true that there is no definitive demonstration that active gravitational mass is
different from inertial mass, there are some experimental anomalies, such as inconsistent measurements of $G$,
that might be interpreted as hints in that direction \cite{Gillies,newG,QLi}.

In this paper, however, it is the equivalence of two relativistic definitions of mass that is called into question instead.
The first is the Einstein mass, defined in the special theory of relativity by
\begin{align}
M:=\sqrt{E^2/c^4-p^2/c^2},
\end{align}
so that it transforms as a scalar under the Lorentz group $SO(3,1)$.
This mass is used in Einstein's theory of gravity (general relativity), so is analogous to, but not necessarily equal to,
the Newtonian active gravitational mass. The second is the Dirac mass, defined by the Dirac equation
\begin{align}
m\psi = \frac{i\hbar}{c}\gamma^\mu\partial_\mu \psi,
\end{align}
so that it transforms as a scalar under the group $SL(2,\CC)$. Since the Dirac mass is used in quantum mechanics,
on which all mechanical forces ultimately depend, this is analogous to, and perhaps equal to, the Newtonian
inertial mass. But it should be noted that the question of equivalence of Einstein and Dirac masses is not the same question as the
equivalence of Newtonian inertial and gravitational masses.

The model proposed in this paper permits the Einstein and Dirac masses to be locally equivalent, but does not permit them to
be globally equivalent. In fact, it is not necessary to invoke this model in order to obtain a mathematical demonstration that these two types of mass
cannot be globally equivalent. This conclusion follows indeed from the general covariance of general relativity, which extends
the Lorentz group $SO(3,1)$ to a group $SL(4,\RR)$. If the Dirac and Einstein masses were equivalent, we would need to
extend the group $SL(2,\CC)$ to a double cover of $SL(4,\RR)$ acting on the Dirac spinor. But no such group exists.

This implies that the local equivalence of Einstein and Dirac masses can only be assumed in a region of spacetime in which
the gravitational field, and/or the non-inertial motion of the experiment, is sufficiently uniform. A number of experiments and
observations in which this condition is not met show significant anomalies, which support the conclusion that
the Einstein and Dirac masses are not the same in these circumstances. In particular, detailed observations of galaxy rotation
curves show that the equivalence of Einstein and Dirac masses breaks down on a galactic scale
\cite{newparadigm}. The flyby anomaly \cite{flyby} may 
similarly be evidence of a breakdown, in circumstances where the relative motion of the satellite and the laboratory is
highly non-inertial. 

If we suppose that the Einstein mass is defined by a copy of $SO(3,1)$ acting on $\rep4a$, representing $4$-momentum,
and we suppose that classical and relativistic physics arises from the quantum effects described by the action of the finite group,
then we must identify the tensors
\begin{align}
\Lambda^2(\rep4a)&=\rep3a+\rep3b,\cr
S^2(\rep4a)&=\rep1+\rep4a+\rep5
\end{align}
for the finite group 
with the corresponding tensors for $SO(3,1)$, that is the representations with spin $(1,0)+(0,1)$ and $(0,0)+(1,1)$ respectively.
Hence the Einstein mass is represented in $\rep1$, and the $5$-momentum (or mass-momentum-energy) lies in $\rep1+\rep4a$,
which is equivalent to the permutation representation of the (binary) icosahedral group on the five cubes inscribed in the
dodecahedron. There is then no consistent action of $SO(3,1)$ on $\rep5$, although one can impose an action that is
equivalent to the action on $\rep1+\rep4a$, if one breaks enough symmetry. The latter action is, however, inconsistent
with general relativity.

If we now also suppose that the Dirac mass is defined by a copy of $SL(2,\CC)$ acting on $\rep4b$, then we can transfer the
action of the compact part $SU(2)$ to an action of $SO(3)$ on $\rep5$ via the standard map from the spin group to the
orthogonal group. We can also transfer the action to the tensors
\begin{align}
\Lambda^2(\rep4b)&=\rep1+\rep5,\cr
S^2(\rep4b)&=\rep3a+\rep3b+\rep4a,
\end{align}
which allows us to extend the action to $SO(5,1)$ on $\rep1+\rep5$, but does not allow us any splitting of
$\rep3a+\rep3b+\rep4a$ into irreducibles, so does not give us any natural action on $\rep3a+\rep3b$.
In particular, the Dirac mass is represented inside $\rep5$, which allows the weak interaction to
change the Dirac mass, while leaving the Einstein mass unchanged.

Hence we see again, more explicitly, that identifying the Dirac mass with the Einstein mass has the effect of
identifying the representations $\rep1+\rep4a$ and $\rep5$. This can only be achieved by breaking the symmetry
of the binary icosahedral group to a subgroup without elements of order $3$. Such subgroups lie either in the
quaternion group $Q_8$ or the dicyclic group of order $20$ (corresponding to the rotations of the icosahedron
that fix an axis between two opposite vertices).

On the other hand, if we allow ourselves to treat the Dirac and Einstein masses separately, then we see the difference
in the difference between the permutation representation $\rep1+\rep4a$ and the monomial representation $\rep5$ of the
icosahedral group. These two representations can be written by mapping the generators $f,g,h$ to the following:
\begin{align}
&(1,2)(3,4),\quad (1,2,3),\quad (1,2)(4,5);\cr
&(\omega1,2)(3,\omega4),\quad (1,\omega2,3)(4,\omega^24,\omega4)(5,\omega5,\omega^25),\quad (1,\omega2)(\omega4,5),
\end{align}
where $\omega$ is a primitive cube root of unity. In particular, writing $\rep5$ as a monomial representation requires
the extension from real to complex numbers. This is perhaps another reason why the standard model requires
a complex rather than real Dirac algebra. However, by identifying these two symmetry groups, we obliterate the
complex numbers in the action of $g$, and thereby obliterate the generation symmetry from the model.

As a final remark, we should note that the definition of mass in the early years of the 20th century was 
clearly the Einstein mass, while the definition of mass in particle physics in the early 21st century is clearly
the Dirac mass. There must therefore have been a changeover period, during which these two distinct types of
mass were calibrated against each other. As I have already indicated, this changeover period appears to have
effectively ended in 1973, and begun probably in the late 1950s. Before this time, mass measurements were
simply not accurate enough for it to be necessary to distinguish Einstein mass from Dirac mass.

\subsection{Implications for quantum gravity}
Attempts to quantise gravity directly, by quantising the gauge group $GL(4,\RR)$, have been unsuccessful \cite{GL4R1,GL4R2}.
The mathematical reason for this is that the finite symmetry groups cannot reach into the non-compact part of the group,
which has the effect of making the theory non-renormalizable. For a Yang--Mills theory of gravity, it is essential for the gauge
group to be compact, although it is not essential for it to be complex unitary, rather than real orthogonal or quaternion
symplectic. In order to obtain a Yang--Mills theory from the discrete model proposed here, it is necessary and
sufficient to take the gauge group
to be
\begin{align}
SO(4) \times SO(5)
\end{align}
acting on $\rep4a+\rep5$,
since $SO(4)$ also acts as $SO(3)\times SO(3)$ on $\rep3a+\rep3b$.

This group is a quotient of the symplectic group
\begin{align}
Sp(1) \times Sp(1) \times Sp(2) & \cong Spin(3) \times Spin(3) \times Spin(5)\cr
& \cong Spin(4)\times Spin(5)
\end{align}
acting on $\rep2a+\rep2b+\rep4b$.
In particular, the mixing of $\rep2a+\rep2b$ with $\rep4b$ that occurs in the standard model
carries with it a mixing of $\rep4a$ with $\rep5$. Note also the similarity between this gauge group and
the Pati--Salam \cite{PatiSalam} gauge group
\begin{align}
SU(2)_L\times SU(2)_R\times SU(4).
\end{align}
The differences are of two kinds: one is the mathematical restriction from the $15$-dimensional $SU(4)$ to the
$10$-dimensional subgroup $Sp(2)$, while the other is a significant difference in physical interpretation.

Now the adjoint representation of our proposed gauge group consists of
\begin{align}
\Lambda^2(\rep4a) & = \rep3a+\rep3b\cr
\Lambda^2(\rep5) & = \rep3a+\rep3b+\rep4a.
\end{align}
Two important things to note here are that, first, the model allows a mixing between the two copies of
$\rep3a+\rep3b$, and second, that $\Lambda^2(\rep5)$ differs from what we would expect from general
relativity, that is
\begin{align}
S^2(\rep4a) = \rep1 + \rep4a+\rep5.
\end{align}
The first property allows us to implement the required mixing between $\rep4a$ and $\rep5$ at the gauge group level,
in order to describe the standard model electroweak mixing of $\rep2a+\rep2b$ and $\rep4b$ at the spinor level.
It also permits a unification of the electromagnetic field in one copy of $\rep3a+\rep3b$
with a $3$-dimensional Newtonian gravitational field and a $3$-dimensional gravito-magnetic field in the other
copy of $\rep3a+\rep3b$. In other words, it permits, and probably requires, a `mixing' between gravity and
electromagnetism at the level of effective quantum field theories.

The second property shows the difference between the permutation and monomial representations of the icosahedral
group on six points. The six points here are the axes of an icosahedron joining opposite vertices, and the difference
between the representations is simply whether one regards a reversed axis as being the same as, or the negative of,
the original axis. 
We may label these axes so that $f,g,h$ act as
\begin{align}
(1,2)(3,4)(5,-5)(6,-6),\quad (1,3,5)(2,-4,6),\quad (3,5)(2,6)(1,-1)(4,-4)
\end{align}
respectively, with or without the signs, as appropriate.
The proposed model therefore inserts these signs, that do not appear in general relativity.

A remarkable consequence of the insertion of these signs is that the hyperspin $2$ representation $\rep5$
disappears from quantum gravity, and with it, apparently, the requirement for a spin $2$ graviton. 
Instead we have the hyperspin $1$ and hyperspin $3$ representations, usually identified as
spin $(1,0)$, spin $(0,1)$ and spin $(1/2,1/2)$.
The first two are normally interpreted, in the context of quantised long-range forces, as spin $1$ photons, while the last
would normally be interpreted as spacetime, or $4$-momentum.
Here we need it to consist of particles, with zero Einstein mass, but possibly non-zero Dirac mass. 
The only reasonable interpretation then seems to be as spin $1/2$ neutrinos and anti-neutrinos.
This permits a quantisation of Euclidean $4$-momentum, and in particular gives the neutrinos
a non-zero Dirac mass, which they need in order to explain neutrino oscillations, together with
a zero Einstein mass, which they need in order to quantise gravity.

However, it is clear that this cannot be the whole story of quantum gravity, because the representation
$\rep5$ contains real physical information that we have not used. In particular, we have taken no
account of neutrino oscillations, which require us to extend from $4$ degrees of freedom,
representing momentum and a non-zero Dirac mass, to $9$, representing momenta for three generations,
combined with a zero Einstein mass. Mathematically, this can be achieved by extending Dirac neutrinos
described by $\rep2a\otimes\rep2b$ to Einstein neutrinos defined by
\begin{align}
\rep3a\otimes\rep3b=\rep4a+\rep5,
\end{align}
but what this means physically is not at all clear.

\section{Symmetry-breaking and the standard model}

\subsection{The standard model}
As I have shown, the standard model of particle physics is based on an identification of two types of Dirac spinors,
in the representations $\rep2a+\rep2b$ and $\rep4b$, that breaks the symmetry of the finite group. There are thus
three different types of `square' of the Dirac spinor, that describe three different ways of looking at the effective fields:
\begin{align}
(\rep2a+\rep2b)\otimes(\rep2a+\rep2b) &= \rep1+\rep4a+\rep3a+\rep3b+\rep4a+\rep1,\cr
(\rep2a+\rep2b)\otimes\rep4b &= \rep5+\rep3a+\rep3b+\rep5,\cr
\rep4b\otimes\rep4b &= \rep1+\rep4a+\rep3a+\rep3b+\rep5.
\end{align}
At scales much greater than the Planck scale, all these representations can be interpreted as (effectively continuous)
fields. But the interpretation may differ according to whether the scale is nuclear, atomic,  classical electromagnetic, or classical
gravitational. All three versions of the square of the spinor contain $\rep3a+\rep3b$, which is usually interpreted as the
electromagnetic field on all scales. 

At the fundamental (Planck scale) level, the last version describes the quantum fields that are predicted by the
model presented in this paper. At the gravitational level, these are essentially the same fields that are predicted by
Einstein's general theory of relativity, although usually interpreted somewhat differently. 
The first version describes the atomic scale, and the standard implementation of QED, while the second
describes the nuclear scale, and the weak and strong forces. Electroweak unification is obtained by mixing
together the first two versions, so that one copy (only) of $\rep1+\rep4a$ is mixed with $\rep5$, via the
identification of
\begin{align}
\rep2a\otimes (\rep2a+\rep2b) &= \rep1+\rep3a+\rep4a,\cr
\rep2a\otimes\rep4b &= \rep3a+\rep5.
\end{align}
This procedure gives the standard model the complete set of effective fields in $\rep4b\otimes \rep4b$ that are needed
in order to give a complete description of physics at the atomic and nuclear scales.

But in order to describe the necessary symmetries on $\rep5$, the standard model interprets $\rep3a+\rep5$ 
and/or $\rep3b+\rep5$ as
the adjoint representation of $SU(3)$, representing the strong force. This procedure unavoidably mixes the strong force
with the electroweak forces which also use $\rep3a$ and/or $\rep3b$.
In other words, the standard model tries to quantise the effective fields separately, but finds that this is impossible,
since the quantised fields are elements of the finite group, that act on all the effective fields simultaneously, whereas the
gauge groups each act on a single field. It is this fundamental distinction between the ways the finite group and
the gauge groups embed in the group algebra that is the mathematical reason why the standard model is forced
to incorporate some very messy `mixing' between the gauge groups, in order to express the relationship
between the gauge fields and the underlying discrete structure.

\subsection{Space and time} 
In order to use this model to describe the dynamics of the universe it is necessary to break the symmetry
of the spacetime representation $\rep4a$ to separate space from time. The obvious way to do this is to
restrict the finite symmetry group to the binary tetrahedral group, generated by $f$ and $g$. But all that is
required mathematically is for $\rep4a$ to become equivalent to either $\rep1+\rep3a$ or $\rep1+\rep3b$,
for which it is necessary and sufficient for the subgroup to contain no elements of order $5$.
Hence the subgroup generated by $g$ and $h$ is another possibility, as is the subgroup generated by $f$ and $h$.
Of course, these three options give three quite different definitions of time, and the last two also distinguish 
one dimension of space from the other two.

Hence only the first possibility could be useful for the purpose of defining time relative to an isotropic $3$-dimensional
physical space. The others may be useful in circumstances where the symmetry of space is broken, for example
by a current, or a gravitational field. In such circumstances, the physical experience of time is known to be different,
and is described by special and general relativity respectively.

Now the group algebra of $\langle f,g\rangle$ is isomorphic to
\begin{align}
\RR + \CC + M_3(\RR) + \HH + M_2(\CC),
\end{align}
with compact subgroup
\begin{align}
U(1) \times SO(3) \times SU(2) \times U(2).
\end{align}
The corresponding real representations have dimensions $1$, $2$, $3$, $4$ and $4$ respectively. The quaternionic $4$, say $4a$, is
the restriction of both $\rep2a$ and $\rep2b$, while $\rep4b$ restricts as two copies of the other $4$, say $4b$.
The quaternionic tensor square of $4a$ is $1+3$, which provides a natural map from $SU(2)$ onto $SO(3)$.
The representation $4b$ has two possible complex structures, related by complex conjugation, so has various different `tensor squares'.
The complex tensor product of the two versions is a complexified copy of $1+3$, that is $1+1+3+3$ as a real representation,
while the complex tensor square of either version is $2+3+3$ as a real representation. This implies that there are \emph{two}
natural maps from $U(2)$ onto $U(1)\times SO(3)$, differing by complex conjugation on the scalars.

In other words, this gives us a toy model with gauge group $SU(2)\times U(2)$, in which the $SU(2)$ factor has a natural interpretation
as the spin group, and $U(2)$ as the gauge group of electroweak theory. However, the natural map from $U(2)$ onto the rotation group
$SO(3)$ is somewhat unexpected here. Taken at face value, it relates the symmetry-breaking of the weak interaction to a symmetry-breaking
of space itself. Whilst this is logically possible, and could, for example, be obtained physically by a mixing between the
weak interaction and quantum gravity, it is not part of the standard model. 

The latter attempts to remove any dependence of the
weak interaction on the gravitational environment by mixing $U(2)$ with $SU(2)$ so that the image of $U(2)$ in $SO(3)$ is cancelled out
by the image of $SU(2)$ in $SO(3)$. Since the direction of spin cannot be measured, this mathematical procedure cannot be detected
experimentally, so that its physical validity can be neither proved nor disproved.
From a mathematical point of view, however, this procedure is invalid, and converts what might have been an exact universal model into
an approximate local model. 

In this approximate local model, however, if we ignore the `mixed' copy of the gauge group, then we effectively identify $SU(2)$ acting on
$4a$ with $SU(2)$ acting on $4b$, so that we can extend this copy of $SU(2)$ to $SL(2,\CC)$ acting on $4b$. Then in order to separate the
left-handed and right-handed spin $1/2$ representations of $SL(2,\CC)$ we must complexify the representation, and hence work
with a complex version of
\begin{align}
4b\otimes 4b = 1+1+3+3+3+3+2
\end{align}
as the Dirac algebra. It is then possible to impose the standard structure of the Dirac algebra as a representation of $SL(2,\CC)$,
but this structure is not determined by the finite group, and is therefore not determined by the fundamental properties of the universe, 
but by our choice of experiments and our choice of models. 

Moreover, by reducing the gauge group from $U(2)$ to $U(1)$, generated in the Feynman calculus by $(1-\gamma_5)i$, we have removed
the possibility of explaining the symmetry-breaking of the weak interaction in terms of the natural map from $U(2)$ onto 
the rotation group $SO(3)$ of space. The finite group model, however, implies that the electroweak mixing angle must be visible as
a parameter of the non-inertial motion of the laboratory, that determines the choice of homomorphism between $SL(2,\CC)$
and $SO(3,1)$ in the standard model. Moreover, there are really three copies of $U(1)$ being mixed together here, so that there
are two distinct physical angles that combine to make the electroweak mixing angle. 
Both are physical angles in space that affect the
quantum gravitational field in which the experiments take place. 

They certainly cannot be angles of latitude or longitude, but must be angles that remain the same for all experiments on the Earth.
They could however be angles that determine tidal effects of gravity, such as the angle of tilt of the Earth's axis, and of
the Moon's orbit, relative to the ecliptic. These angles are approximately $23.44^\circ$ and $5.14^\circ$ respectively,
and their sum, $28.58^\circ$ is indeed very close to the various different values of the electroweak mixing angle
that are measured in different experiments. This coincidence is not exactly a prediction of the model,
since there are many things I have not taken into account in this toy model. Nevertheless, it is striking,
and suggests that it might be fruitful to examine this model more closely.

\section{Conclusion}
In this paper I have 
argued that there is essentially only one algebraic structure in which it is 
conceivable to build
a mathematical model of quantum mechanics that 
is consistent with the most fundamental physical principles that are required by experiment, including
locality, unitarity, gauge invariance and renormalisability.
The structure of quaternion reflections
leads naturally to both 
first and second quantisation, 
which appear to go beyond the standard model only by tripling the numbers of particle states
for electrons and photons, and reinterpreting the colourless gluons as quantum superpositions
of the `primordial' massless weak bosons. 
Apart from this, no new particles or new forces are required.

I have built only the basic algebraic structures,
and have not put the physics on top.
There is therefore no guarantee that a model of this type can include the standard model,
and no guarantee that it will be useful. 
Nevertheless it contains some important ingredients that will be necessary 
in any  theory that seeks to go beyond the standard model in a meaningful way, including 
\begin{itemize}
\item a mathematical definition of the elementary fermions that includes all 
three generations, 
\item a mathematical process whereby the elementary fermions generate quantum fields, 
\item a mathematical definition of massless mediators that includes photons and gluons, 
\item a mathematical definition of massive mediators that includes all pairs of charged
intermediate vector bosons and pseudoscalar mesons, and
\item an embedding of the gauge groups of the standard model into larger gauge groups that take account of the
three generations. 
\end{itemize}
In addition, it contains
\begin{itemize}
\item
enough quantum states to explain the experimental measurement of electron spin and photon polarisation
without the need for superluminal transfer of information, and
\item a mathematical isomorphism between the quantum and 
classical fields. 
\end{itemize}
This last 
suggests a division of both into dimensions $3+6+6$, identified as the God-given fundamental
gravitational, electromagnetic
and colour fields respectively. 
But in passing from the fundamental theory to the effective field theories we pass from the structure of the group itself
to the structure of the group algebra, converted into a Lie algebra, which divides the fields instead as $3+3+4+5$.
The relationship between the group \emph{elements} that parametrise $3+6+6$ and the group \emph{representations}
that parametrise $3+3+4+5$ is quite complicated.

In particular, to see the relationship  it is necessary to break the symmetry and 
mix the effective fields $3+3+4+5$ together in complicated ways in order to obtain
the fundamental quantum fields in $3+6+6$. All this symmetry-breaking and mixing has nothing
whatever to do with the fundamental theory, and has only to do with our choice of effective field theories.
Mainstream opinion is that our effective field theories are fundamental, but
by examining the data, I provide evidence that they may not be, and that we may have a lot more choice than is generally realised,
for parameters to use 
in constructing effective field theories. Chief amongst these choices is our choice of the concept of
mass, specifically our choice to identify the Einstein mass with the Dirac mass.

Finally, let me remark that this work does not solve the mathematical Yang--Mills mass gap problem as described in detail
in \cite{Clay}. It does however potentially solve the underlying physical problem, by modelling the discrete nature of mass. As such, it
may
render any potential solution of the Clay Millennium problem largely irrelevant for the understanding of the
physical universe.


\begin{thebibliography}{99}

\bibitem{YangMills} C. N. Yang and R. L. Mills (1954),
Conservation of isotopic spin and isotopic gauge invariance,
{\it Phys. Rev.} {\bf 96}, 191.

\bibitem{tHooft22} G. 't Hooft (2022), Projecting local and global symmetries to the Planck scale,
arXiv:2202.05367.

\bibitem{Blichfeldt} H. F. Blichfeldt (1917), {\it Finite collineation groups}, The University of Chicago Press.


\bibitem{finite} R. A. Wilson (2021), Finite symmetry groups in physics, arXiv:2102.02817.

\bibitem{octahedral} R. A. Wilson (2021), Options for a finite group model of quantum mechanics,
arXiv:2104.10165.

\bibitem{gl23model}
R. A. Wilson (2021), 
Potential applications of modular representation theory to quantum mechanics,
arXiv:2106.00550.

\bibitem{icosa} R. A. Wilson (2021),
Possible uses of the binary icosahedral group in grand unified theories, arXiv:2109.06626.



\bibitem{Porteous} I. R. Porteous (1995), {\it Clifford algebras and the classical groups},
Cambridge UP.

\bibitem{Cohen} A. M. Cohen (1980), Finite quaternionic reflection groups, 
{\it J. Algebra} {\bf 64}, 293--324.

\bibitem{icosians} R. A. Wilson (1986), The geometry of the Hall--Janko group as a quaternionic reflection group,
{\it Geom. Dedicata} {\bf 20}, 157--173.

\bibitem{Bell} J. S. Bell (1964), On the Einstein Podolsky Rosen paradox, {\it Physics} {\bf 1}, 195--200.


\bibitem{Rae} A. Rae (1986), {\it Quantum physics: illusion or reality?}, Cambridge UP.


\bibitem{WoitQFT} P. Woit (2017), {\it Quantum theory, groups and representations}, Springer.



 \bibitem{Zee}
 A. Zee (2016), {\it Group theory in a nutshell for physicists},
 Princeton University Press.
  

\bibitem{JamesLiebeck} G. D. James and M. W. Liebeck (2012), {\it Representations and characters of groups}, 2nd ed,
Cambridge UP.


\bibitem{Klein} F. Klein (1884), {\it Lectures on the ikosahedron}, English translation, 1888.


\bibitem{Klein1} F. Klein (1876), {\it Math. Annalen} {\bf 9}, 183.

\bibitem{DuVal} P. Du Val (1964), {\it Homographies, quaternions and rotations}, Clarendon Press, Oxford.


\bibitem{SPLAG}
J. H. Conway and N. J. A. Sloane (1988), {\it Sphere packings, lattices and groups}, Springer. 

\bibitem{tHooft17} G. 't Hooft (2017), Free will in the theory of everything, arXiv:1709.02874.


\bibitem{almanac}  {\it The astronomical almanac}, and supplement (1961).


\bibitem{1973} E. R. Cohen and B. N. Taylor (1973),
The 1973 least-squares adjustments of the fundamental constants,
{\it J. Phys. Chem. Ref. Data} {\bf 2} (4), 663--734.

\bibitem{Gillies} G. T. Gillies (1997), The Newtonian gravitational constant: recent measurements and related studies,
{\it Reports on progress in physics} {\bf 60} (2), 151--225.

\bibitem{newG} C. Rothleitner and S. Schlamminger (2017),
Invited review article: Measurements of the Newtonian constant of gravitation, $G$,
{\it Review of Scientific Instruments} {\bf 88}, 111101.

\bibitem{QLi} Q. Li et al. (2018), Measurements of the gravitational constant using two independent methods,
{\it Nature} {\bf 560} (7720): 582--588. 

\bibitem{newparadigm} P. Kroupa, M. Pawlowski and M. Milgrom (2012),
The failures of the standard model of cosmology require a new paradigm,
{\it Intern. J. Modern Phys. D} {\bf 21} (14).
arXiv:1301.3907.


\bibitem{flyby} J. D. Anderson, J. K. Campbell, J. E. Ekelund, J. Ellis and J. F. Jordan (2008),
Anomalous orbital-energy changes observed during spacecraft flybys of Earth,
{\it PRL} {\bf 100}, 091102.


\bibitem{GL4R1}
D. Ivanenko and G. Sardanashvily (1983), The gauge treatment of gravity, {\it Physics Reports} {\bf 94}, 1--45.

\bibitem{GL4R2} M. Leclerc (2006), The Higgs sector of gravitational gauge theories, {\it Annals of Physics} {\bf 321}, 708--743.

\bibitem{PatiSalam} J. C. Pati and A. Salam (1974), Lepton number as the fourth `color',
{\it Phys. Rev. D} {\bf 10} (1), 275--289.

\bibitem{Clay} A. Jaffe and E. Witten, Quantum Yang--Mills theory,
www.claymath.org/sites/default/files/yangmills.pdf

\end{thebibliography}
\end{document}